# Variable selection in measurement error models

YANYUAN MA[1] and RUNZE LI[2]

[1]*Department of Statistics, Texas A&M University, College Station, TX 77843, USA.*
*E-mail: ma@stat.tamu.edu*

[2]*Department of Statistics and The Methodology Center, The Pennsylvania State University, University Park, PA 16802, USA. E-mail: rli@stat.psu.edu*

Measurement error data or errors-in-variable data have been collected in many studies. Natural criterion functions are often unavailable for general functional measurement error models due to the lack of information on the distribution of the unobservable covariates. Typically, the parameter estimation is via solving estimating equations. In addition, the construction of such estimating equations routinely requires solving integral equations, hence the computation is often much more intensive compared with ordinary regression models. Because of these difficulties, traditional best subset variable selection procedures are not applicable, and in the measurement error model context, variable selection remains an unsolved issue. In this paper, we develop a framework for variable selection in measurement error models via penalized estimating equations. We first propose a class of selection procedures for general parametric measurement error models and for general semi-parametric measurement error models, and study the asymptotic properties of the proposed procedures. Then, under certain regularity conditions and with a properly chosen regularization parameter, we demonstrate that the proposed procedure performs as well as an oracle procedure. We assess the finite sample performance via Monte Carlo simulation studies and illustrate the proposed methodology through the empirical analysis of a familiar data set.

## 1. Introduction

In the regression analysis, some covariates often can only be measured imprecisely or indirectly, thus resulting in measurement error models, also known as errors-in-variable models in the literature. Various statistical procedures have been developed for statistical inference in measurement error models (Carroll, Ruppert, Stefanski and Crainiceanu (2006)). The study on linear measurement error models dates back to Bickel and Ritov (1987), where an efficient estimator is provided. Stefanski and Carroll (1987) constructed consistent estimators for generalized linear measurement error models. Recently, Tsiatis and Ma (2004) extended the model framework to an arbitrary parametric regression setting. Liang, Härdle and Carroll (1999) proposed partially linear measurement error models. Ma and Carroll (2006) studied generalized partially linear measurement error models. Further active research development has been established recently in the nonparametric measurement error area; see, for example, Delaigle and Hall (2007) and Delaigle and Meister (2007). The goal of this paper is to develop a class of variable selection procedures for general measurement error models. We would emphasize here that the scope of the paper is not limited to generalized linear models.

This study was motivated by examining the effects of systolic blood pressure (SBP), a covariate with error, and the effects of three other covariates – respectively, serum cholesterol, age, and smoking status – on the probability of the occurrence of heart disease. In our initial analysis, we include interactions between SBP and covariates and interactions among covariates and quadratic terms of covariates to reduce modeling bias. It is found in our preliminary analysis that some interactions and quadratic terms are not significant and should be excluded to achieve a parsimonious model. To select significant variables in further analysis, we realized that the traditional Akaike information criterion (AIC) and Bayesian information criterion (BIC) criteria are not well defined for the model we consider in Section 4.4. Recently, a class of variable selection procedures for partially linear measurement error models via using penalized least squares and penalized quantile regression were proposed in Liang and Li (2009). However, their procedures are not applicable to cases beyond partially linear models, such as partially linear logistic regression models, and therefore the procedures in Liang and Li (2009) cannot be applied for the model in Section 4.4 either. In fact, variable selection for general parametric or semiparametric measurement error models is challenging. One major difficulty is the lack of a likelihood function in these models, due to the difficulty in obtaining the distribution of the error-prone covariates. For example, using $Y$ to denote the response variable, $X$ to denote the unobservable covariate, and $W$ to denote an observed surrogate of $X$, the likelihood of a single observation $(w, y)$ is then $\int p_{Y|X}(y|x, \beta) p_{W|X}(w|x) p_X(x) \, \mathrm{d}x$. In order to calculate this likelihood, one will need to estimate $p_X$, yielding a deconvolution problem that is known to have a very slow rate (Carroll and Hall (1988), Fan (1991)) and is typically avoided in parametric measurement error models. Although a reasonable criterion function can be used in place of the likelihood, the difficulty persists in that, except for very special models such as in linear or partially linear cases, even a sensible criterion function is unavailable. In other models such as the ones that arise in survival analysis, the lack of a likelihood function also causes a problem. To perform variable selection in these models, rather complicated methods have been proposed where for each potential model, one needs to fit the model, derive the asymptotic properties of the estimator, form some artificial criterion function based on the asymptotic properties of the estimators, and finally add a penalty to perform the procedure. The procedure is complicated and unnatural. This motivates us to develop some simple variable selection procedures for measurement error models when a reasonable criterion function is unavailable. Although a few variable selection procedures exist for linear or partially linear measurement error models (Liang and Li (2009)), to the best of our knowledge, variable selection for general parametric or semi-parametric measurement error models has never been systematically studied in the literature. This paper intends to fill this gap by developing a class of variable selection procedures for both parametric and semi-parametric measurement error models. In addition, the method proposed here is applicable to the more general situation where the likelihood or a natural criterion is not available, and the estimation is performed through solving a set of estimating equations.

The variable selection procedure we propose is indeed a penalized estimating equation method that can be applied for both parametric and semi-parametric measurement error models. In addition, the penalized estimating equation method is applicable to any set of consistent estimating equations. Note that here the measurement error model we consider is completely general and not limited to generalized linear models. Variable selection and feature selection are very active research topics in the current literature. Candès and Tao (2007) and Fan and Lv (2008) have studied variable selection for linear models when the sample size is much smaller than the dimension of the regression parameter space. Their results are inspiring, but only valid for linear models with very strong assumptions on the design matrix or the distribution of covariates. Thus, in this paper we follow Fan and Peng (2004) and consider the statistical setting in which the number of regression coefficients diverges to infinity at a certain rate as the sample size tends to infinity. We systematically study the asymptotic properties of the proposed estimator. It is worth pointing out that theoretic results in this paper provide explicit results on the asymptotic properties when the dimension of regression coefficients increases as the sample size increases. This advances the results in current literature, where estimation and inference are studied only for fixed finite-dimensional parameters for measurement error models. In our asymptotic analysis, we show that with a proper choice of the regularization parameters and the penalty function, our estimator possesses the oracle property, which roughly means that the estimate is as good as when the true model is known (Fan and Li (2001)). We also demonstrate that the oracle property holds in a simpler form for the more familiar setting where the true number of regression coefficients is fixed.

In addition, we address issues of practical implementation of the proposed methodology. It is desirable to have an automatic, data-driven method to select the regularization parameters. To this end, we propose generalized cross-validation (GCV)-type and BIC-type tuning parameter selectors for the proposed penalized estimating equation method. Monte Carlo simulation studies are conducted to assess finite sample performance in terms of model complexity and model error. From our simulation studies, both tuning parameter selectors result in sparse models, while the BIC-type tuning parameter selector outperforms the GCV-type tuning parameter selectors.

The rest of the paper is organized as follows. In Section 2, we propose a new class of variable selection procedures for parametric measurement error models and study asymptotic properties of the proposed procedures. We develop a new variable selection procedure for semi-parametric measurement error models in Section 3. Implementation issues and numerical examples are presented in Section 4, where we describe data-driven automatic tuning parameter selection methods (Section 4.1), define the concept of approximate model error to evaluate the selected model (Section 4.2), carry out a simulation study to assess the finite sample performance of the proposed procedures (Section 4.3), and illustrate our method in an example (Section 4.4). Technical details are collected in the Appendix.

## 2. Parametric measurement error models

A general parametric measurement error model has two parts, written as

$$p_{Y|X,Z}(Y|X,Z,\beta) \quad \text{and} \quad p_{W|X,Z}(W|X,Z,\xi). \tag{2.1}$$

Here, the main model is $p_{Y|X,Z}(Y|X,Z,\beta)$, which denotes the conditional probability density function (p.d.f.) of the response variable $Y$ on the covariates measured with error $X$ and the covariates measured without error $Z$. Note that here the conditional distribution of $Y$ on the covariates is completely general, hence it includes many familiar regression families. For example, the linear model with normal error $Y = X\beta + e$, the logistic model $\Pr(Y=1|X) = \exp(\beta_0 + X\beta_1)/\{1 + \exp(\beta_0 + X\beta_1)\}$, or the Poisson model $p_{Y|X}(Y|X) = \exp\{-(\beta_0 + X^{\mathrm{T}}\beta_1)\}(\beta_0 + X^{\mathrm{T}}\beta_1)^Y/Y!$ are all special forms of the model we consider. The error model is denoted $p_{W|X,Z}(W|X,Z,\xi)$, where $W$ is an observable surrogate of $X$. The parameter $\beta$ is a $d$-dimensional regression coefficient, $\xi$ is a finite-dimensional parameter, and our main interest is in selecting the relevant subset of covariates in $X$ and $Z$ and estimating the subsequent parameters contained in $\beta$. Typically, $\xi$ is a nuisance parameter and its estimation usually requires multiple observations or instrumental variables. As routinely done in the literature, we assume that the model is identifiable. Furthermore, for simplicity, we assume in the main context of this paper that the error model $p_{W|X,Z}(W|X,Z)$ is completely known and hence $\xi$ is suppressed. The extension to the unknown $\xi$ case is rather straightforward and is discussed in Section 5. The observed data is of the form $\{(W_i, Z_i, Y_i), i = 1, \ldots, n\}$.

Denote $S_\beta^*$ as the purported score function. That is,

$$S_\beta^*(W,Z,Y) = \frac{\partial \log \int p_{W|X,Z}(W|X,Z) p_{Y|X,Z}(Y|X,Z) p_{X|Z}^*(X|Z)\,\mathrm{d}\mu(X)}{\partial \beta},$$

where $p_{X|Z}^*(X|Z)$ is a conditional p.d.f. that one posits, which can be either equal or not equal to the true p.d.f. $p_{X|Z}(X|Z)$. Let the function $a(X,Z)$ satisfy

$$E[E^*\{a(X,Z)|W,Z,Y\}|X,Z] = E\{S_\beta^*(W,Z,Y)|X,Z\},$$

where $E^*$ indicates that the expectation is calculated using the posited $p_{X|Z}^*(X|Z)$. Note that here and in the sequel a model $p_{X|Z}^*(X|Z)$ has to be proposed in order to actually construct the estimator. Define

$$S_{\mathit{eff}}^*(W,Z,Y) = S_\beta^*(W,Z,Y) - \mathrm{E}^*\{a(X,Z)|W,Z,Y\}.$$

To select significant variables and estimate the corresponding parameters simultaneously, we propose the penalized estimating equations for model (2.1) as

$$\sum_{i=1}^{n} S_{eff}^{*}(W_i, Z_i, Y_i, \beta) - n\dot{p}_\lambda(\beta) = 0, \tag{2.2}$$

where $\dot{p}_\lambda(\beta) = \{p'_\lambda(\beta_1), \ldots, p'_\lambda(\beta_d)\}^{\mathrm{T}}$ and $p'_\lambda(\cdot)$ is the first-order derivative of a penalty function $p_\lambda(\cdot)$. Solving for $\hat{\beta}$ from (2.2) gives the estimate of $\beta$. In practice, we may allow different coefficients to have penalty functions with different regularization parameters. For example, we may want to keep some variables in the model without penalizing their coefficients. For ease of presentation, we assume that the penalty functions and the regularization parameters are the same for all coefficients in this paper.

The penalties in the classic variable selection criteria, such as AIC and BIC, cannot be applied to the penalized estimating equations. Following the study on the choice of the penalty functions in Fan and Li (2001), we use the smoothly clipped absolute deviance (SCAD) penalty, whose first-order derivative is defined as

$$p'_\lambda(\gamma) = \lambda\left\{I(|\gamma| \leq \lambda) + \frac{(a\lambda - |\gamma|)_+}{(a-1)\lambda}I(|\gamma| > \lambda)\right\}\mathrm{sign}(\gamma) \tag{2.3}$$

for any scalar $\gamma$, where $\mathrm{sign}(\cdot)$ is the sign function, that is, $\mathrm{sign}(\gamma) = -1$, 0 and 1 when $\gamma < 0$, $= 0$ and $> 0$, respectively. Here, $a > 2$ is a constant, and a choice of $a = 3.7$ is appropriate from a Bayesian point of view. A property of (2.2) is that with a proper choice of penalty functions, such as the SCAD penalty, the resulting estimate contains some exact zero coefficients. This is equivalent to excluding the corresponding variables from the final selected model, thus variable selection is achieved at the same time as parameter estimation.

Concerns about model bias often prompt us to build models that contain many variables, especially when the sample size becomes large. A reasonable way to capture such a tendency is to consider the situation where the dimension of the parameter $\beta$ increases along with the sample size $n$. We therefore study the asymptotic properties of the penalized estimating equation estimator under the setting in which both the dimension of the true non-zero components of $\beta$ and the total length of $\beta$ tend to infinity as $n$ goes to infinity. Denote $\beta_0 = (\beta_{10}, \ldots, \beta_{d0})^{\mathrm{T}}$ as the true value of $\beta$. Let

$$a_n = \max\{|p'_{\lambda_n}(|\beta_{j0}|)| : \beta_{j0} \neq 0\} \quad \text{and} \quad b_n = \max\{|p''_{\lambda_n}(|\beta_{j0}|)| : \beta_{j0} \neq 0\}, \tag{2.4}$$

where we write $\lambda$ as $\lambda_n$ to emphasize its dependence on the sample size $n$.

**Theorem 1.** *Suppose that condition (*P1*) in the Appendix holds. Under regularity conditions (*A1*)–(*A3*) in the Appendix, and if $d_n^4 n^{-1} \to 0$, $\lambda_n \to 0$ when $n \to \infty$, then with probability tending to one, there exists a root of (2.2), denoted $\hat{\beta}$, such that $\|\hat{\beta} - \beta_0\| = \mathrm{O}_p\{\sqrt{d_n}(n^{-1/2} + a_n)\}$, where we write $d$ as $d_n$ to emphasize its dependence on the sample size $n$.*

The proof of Theorem 1 is given in the Appendix. Theorem 1 demonstrates that the convergence rate depends on the penalty function and the regularization parameter $\lambda_n$ through $a_n$. From Theorem 1, it requires $a_n = \mathrm{O}(1/\sqrt{n})$ to achieve root $(n/d_n)$ convergence rate. For the $L_1$ penalty, $a_n = \lambda_n$. Thus, the root $(n/d_n)$ convergence rate requires that $\lambda_n = \mathrm{O}(1/\sqrt{n})$, while $a_n = 0$ as $\lambda_n \to 0$ for the SCAD penalty. Thus, the resulting estimate with the SCAD penalty is root $(n/d_n)$ consistent.

To present the oracle property of the resulting estimate, we first introduce some notation. Without loss of generality, we assume $\beta_0 = (\beta_{I0}^{\mathrm{T}}, \beta_{II0}^{\mathrm{T}})^{\mathrm{T}}$, and in the true model any element in $\beta_{I0}$ is not equal to 0 while $\beta_{II0} \equiv 0$. Denote the dimension of $\beta_I$ as $d_1$. Furthermore, denote

$$b = \{p'_{\lambda_n}(\beta_{10}), \ldots, p'_{\lambda_n}(\beta_{d_10})\}^{\mathrm{T}} \quad \text{and} \quad \Sigma = \mathrm{diag}\{p''_{\lambda_n}(\beta_{10}), \ldots, p''_{\lambda_n}(\beta_{d_10})\}, \qquad (2.5)$$

and the first $d_1$ components of $S^*_{eff}(W, Z, Y, \beta)$ as $S^*_{eff,I}(\beta)$. In the following theorem, we use the same formulation as that in Cai, Fan, Li and Zhou (2005).

**Theorem 2.** *Suppose that condition (*P1*) holds. Under regularity conditions (*A1*)–(*A3*), assume $\lambda_n \to 0$ and $d_n^5/n \to 0$ when $n \to \infty$. If*

$$\liminf_{n\to\infty} \liminf_{\gamma\to 0+} \sqrt{n/d_n}\, p'_{\lambda_n}(\gamma) \to \infty, \qquad (2.6)$$

*then with probability tending to one, any root $n$ consistent solution $\hat{\beta}_n = (\hat{\beta}_I^{\mathrm{T}}, \hat{\beta}_{II}^{\mathrm{T}})^{\mathrm{T}}$ of (2.2) must satisfy that:*

(i) $\hat{\beta}_{II} = 0$,
(ii) *for any $d_1 \times 1$ vector $v$, s.t. $v^{\mathrm{T}} v = 1$,*

$$\sqrt{n} v^{\mathrm{T}} [E\{S^*_{eff,I}(\beta_{I0}) S^{*\mathrm{T}}_{eff,I}(\beta_{I0})\}]^{-1/2} \left\{ E \frac{\partial S^*_{eff,I}(\beta_{I0})}{\partial \beta_I^{\mathrm{T}}} - \Sigma \right\}$$

$$\times \left[ \hat{\beta}_I - \beta_{I0} - \left\{ E \frac{\partial S^*_{eff,I}(\beta_{I0})}{\partial \beta_I^{\mathrm{T}}} - \Sigma \right\}^{-1} b \right] \xrightarrow{\mathrm{D}} N(0,1),$$

*where the notation $\xrightarrow{\mathrm{D}}$ stands for convergence in distribution.*

The proof of Theorem 2 is given in the Appendix. For some penalty functions, including the SCAD penalty, $b$ and $\Sigma$ are zero when $\lambda_n$ is sufficiently small, hence the results in Theorem 2 imply that the proposed procedure has the celebrated oracle property: that is, $\hat{\beta}_{II} = 0$, and for any $d_1 \times 1$ vector $v$, s.t. $v^{\mathrm{T}} v = 1$,

$$\sqrt{n} v^{\mathrm{T}} [E\{S^*_{eff,I}(\beta_{I0}) S^{*\mathrm{T}}_{eff,I}(\beta_{I0})\}]^{-1/2} E\left\{ \frac{\partial S^*_{eff,I}(\beta_{I0})}{\partial \beta_I^{\mathrm{T}}} \right\} (\hat{\beta}_I - \beta_{I0}) \xrightarrow{\mathrm{D}} N(0,1). \quad (2.7)$$

Theorems 1 and 2 imply that for fixed and finite $d$, $\|\hat{\beta} - \beta_0\| = \mathrm{O}_p(n^{-1/2} + a_n)$ and with probability tending to one, any root $n$ convergence solution $\hat{\beta} = (\hat{\beta}_I^{\mathrm{T}}, \hat{\beta}_{II}^{\mathrm{T}})^{\mathrm{T}}$ of (2.2) must satisfy that $\hat{\beta}_{II} = 0$ and

$$\sqrt{n}\bigg[\hat{\beta}_I - \beta_{I0} - \bigg\{E\frac{\partial S^*_{eff,I}(\beta_{I0})}{\partial \beta_I^{\mathrm{T}}} - \Sigma\bigg\}^{-1} b\bigg]$$
$$\xrightarrow{\mathrm{D}} N\bigg[0, \bigg\{E\frac{\partial S^*_{eff,I}(\beta_{I0})}{\partial \beta_I^{\mathrm{T}}} - \Sigma\bigg\}^{-1} E\{S^*_{eff,I}(\beta_{I0})S^{*\mathrm{T}}_{eff,I}(\beta_{I0})\}$$
$$\times \bigg\{E\frac{\partial S^*_{eff,I}(\beta_{I0})}{\partial \beta_I^{\mathrm{T}}} - \Sigma\bigg\}^{-\mathrm{T}}\bigg],$$

where the notation $M^{-\mathrm{T}}$ denotes $(M^{-1})^{\mathrm{T}}$ for a matrix $M$. These results are still valid under much weaker conditions. See an elaborated version of this paper, Ma and Li (2007), for details.

For SCAD penalty and for fixed and finite $d$, (2.7) becomes

$$\sqrt{n}(\hat{\beta}_I - \beta_{I0}) \to N\{0, E(\partial S^*_{eff,I}/\partial \beta_I^{\mathrm{T}})^{-1} E(S^*_{eff,I} S^{*\mathrm{T}}_{eff,I}) E(S^*_{eff,I}/\partial \beta_I^{\mathrm{T}})^{-\mathrm{T}}\}$$

in distribution. In other words, with probability tending to 1, the penalized estimator performs in the same manner as the locally efficient estimator under the correct model.

## 3. Semi-parametric measurement error models

To motivate the problems considered in this section, we start with some commonly used semi-parametric regression models for which the proposed procedure in this section can be directly applied. Consider first the error-free regression cases, and let $Y$ be the response and $Z$ and $S$ be the covariates. Throughout this paper, we consider univariate $Z$ only. Consider the partially linear model defined as follows:

$$Y = \theta(Z) + S^{\mathrm{T}}\beta + \varepsilon. \tag{3.1}$$

The partially linear model keeps the flexibility of the nonparametric model for the baseline function while maintaining the explanatory power of parametric models. Therefore, it has received a lot of attention in the literature. See, for example, Härdle, Liang and Gao (2000) and references therein. Various extensions of the partially linear model have been proposed in the literature. Li and Nie (2007, 2008) proposed the partially nonlinear models

$$Y = \theta(Z) + f(S;\beta) + \varepsilon, \tag{3.2}$$

where $f(S;\beta)$ is a specific, known function that may be nonlinear in $\beta$. See Li and Nie (2007, 2008) for some interesting examples. Li and Liang (2008) and Lam and Fan (2008) studied the generalized varying coefficient partially linear model

$$g\{E(Y|Z,S)\} = S_1^{\mathrm{T}}\beta + S_2^{\mathrm{T}}\theta(Z), \tag{3.3}$$

where $g(\cdot)$ is a link function and $(S_1, S_2, Z)$ are covariates. Model (3.3) includes most commonly used semi-parametric models, such as the partially linear models (3.1), the generalized partially linear models (Severini and Staniswalis (1994)), and semi-varying coefficient partially linear models (Fan and Huang (2005)).

In the presence of covariates measured with error, one may extend the aforementioned semi-parametric regression models for measurement error data. As in the last section, let $X$ be the covariate vector measured with error. Among these semi-parametric models with error, the partially linear measurement error model

$$Y = \theta(Z) + X^{\mathrm{T}}\beta_1 + S^{\mathrm{T}}\beta_2 + \varepsilon \qquad (3.4)$$

has been studied in Liang, Härdle and Carroll (1999). Liang and Li (2009) proposed a class of variable selection for model (3.4) using penalized least squares and penalized quantile regression. Our procedure in this section, however, is directly applicable for both the generalized varying coefficient partially linear measurement error model

$$g\{E(Y|X,Z,S)\} = X^{\mathrm{T}}\beta_1 + S_1^{\mathrm{T}}\beta_2 + S_2^{\mathrm{T}}\theta(Z), \qquad (3.5)$$

where $S = (S_1^{\mathrm{T}}, S_2^{\mathrm{T}})^{\mathrm{T}}$, and the partially nonlinear measurement error model

$$Y = \theta(Z) + f(X, S; \beta) + \varepsilon, \qquad (3.6)$$

when an error distribution is assumed. It is worth noting that model (3.6) includes the following model as a special case

$$Y = X^{\mathrm{T}}\beta_1 + S^{\mathrm{T}}\beta_2 + (XZ)^{\mathrm{T}}\beta_3 + \theta(Z) + \varepsilon, \qquad (3.7)$$

where $(XZ)$ consists of all interaction terms between $X$ and $Z$, but model (3.4) does not. Thus, the variable selection procedures proposed in Liang and Li (2009) are not directly applicable for model (3.7).

In summary, in this section, we consider a general semi-parametric error model that includes models (3.4)–(3.6) as its special cases. Specifically, the semi-parametric measurement error model we consider here also has two parts:

$$p_{Y|X,Z,S}\{Y|X,Z,S,\beta,\theta(Z)\} \quad \text{and} \quad p_{W|X,Z,S}(W|X,Z,S). \qquad (3.8)$$

The major difference from its parametric counterpart is that the main model contains unknown functions $\theta(Z)$. It is easy to check that models (3.4)–(3.6) are special cases of model (3.8). Note that a simpler version of this model is considered in Ma and Carroll (2006), where the dimension of $\beta$ is assumed to be fixed and the dimension of $\theta$ is assumed to be one.

Throughout this paper, we assume that the model is identifiable. We propose the penalized estimating equation for the semi-parametric model:

$$\sum_{i=1}^{n} \mathcal{L}(W_i, Z_i, S_i, Y_i, \beta, \hat{\theta}_i) - n\dot{p}_\lambda(\beta) = 0, \qquad (3.9)$$

where $\dot{p}_\lambda(\beta)$ has the same form as in (2.3). However, the computation of $\mathcal{L}$ is more involved. Denote the dimension of $\theta(Z)$ as $m$, a fixed and finite integer. If we replace $\theta(Z)$ with a single unknown $m$-dimensional parameter $\alpha$ and append $\alpha$ to $\beta$, we obtain from (3.8) a parametric measurement error model with parameters $(\beta^{\mathrm{T}}, \alpha^{\mathrm{T}})^{\mathrm{T}}$. For this parametric model, we can compute the corresponding $S^*_{\mathit{eff}}$ as done in the last section. Specifically, we will have $S^*_{\mathit{eff}}(W, Z, S, Y) = S^*_{\beta,\alpha}(W, Z, S, Y) - E^*\{a(X, Z, S)|W, Z, S, Y\}$, where $a(X, Z, S)$ satisfies $E[E^*\{a(X, Z, S)|W, Z, S, Y\}|X, Z, S] = E\{S^*_{\beta,\alpha}(W, Z, S, Y)|X, Z, S\}$. Note that $S^*_{\mathit{eff}}$ has the same dimension as the dimension of $\beta$ plus $m$. We write the last $m$ components of $S^*_{\mathit{eff}}$ as $\Psi(X, Z, S, Y, \beta, \alpha)$ and the rest as $\mathcal{L}(X, Z, S, Y, \beta, \alpha)$. We now solve for $\hat{\theta}_i$, $i = 1, \ldots, n$, from

$$\sum_{i=1}^n K_h(z_i - z_1)\Psi(w_i, z_i, s_i, y_i; \beta, \theta_1) = 0,$$
$$\vdots \qquad (3.10)$$
$$\sum_{i=1}^n K_h(z_i - z_n)\Psi(w_i, z_i, s_i, y_i; \beta, \theta_n) = 0,$$

where $K_h(z) = h^{-1}K(z/h)$, $K$ is a smooth symmetric kernel function with compact support that satisfies $\int K(t)t^2 \, dt = 1$, and $h$ is a bandwidth. Note that $\theta_1, \ldots, \theta_n$ are all $m$-dimensional parameters. Inserting the $\hat{\theta}_i$'s into $\mathcal{L}$ in (3.9), we obtain a complete description of the estimator. Note that $\hat{\theta}_i$ depends on $\beta$, so a more precise notation for $\hat{\theta}_i$ is $\hat{\theta}_i(\beta)$. Solving equation (3.9) yields a penalized estimating equation estimate. Theorem 3 below gives its convergence rate.

**Theorem 3.** *Suppose that condition (*P1*) holds. Under regularity conditions (*B1*)–(*B4*) in the Appendix, and if $d_n^4 n^{-1} \to 0$, $\lambda_n \to 0$ when $n \to \infty$, then with probability tending to one, there exists a root of (3.9), denoted $\hat{\beta}$, such that $\|\hat{\beta} - \beta_0\| = \mathrm{O}_p\{\sqrt{d_n}(n^{-1/2} + a_n)\}$.*

The proof of Theorem 3 is given in the Appendix. Theorem 3 indicates that to achieve root $(n/d_n)$ convergence rate (or root $n$ convergence rate for finite and fixed $d$), $\lambda_n$ and the penalty function must be chosen such that $a_n = \mathrm{O}_p(n^{-1/2})$.

Let $\mathcal{L}_I$ be the first $d_1$ components of $\mathcal{L}$, $\mathcal{L}_{I\beta_I}$ the partial derivative of $\mathcal{L}_I$ with respect to $\beta_I$, $\mathcal{L}_{I\theta}$ the partial derivative of $\mathcal{L}_I$ with respect to $\theta$, $\Psi_\theta$ the partial derivative of $\Psi$ with respect to $\theta$, and $\Psi_{\beta_I}$ the partial derivative of $\Psi$ with respect to $\beta_I$. Also define $\Omega(Z) = E(\Psi_\theta|Z)$, $\mathcal{U}_I(Z) = E(\mathcal{L}_{I\theta}|Z)\Omega^{-1}(Z)$ and $\theta_{\beta_I}(Z) = -\Omega^{-1}(Z)E(\Psi_{\beta_I}|Z)$. Further defining

$$A = E[\mathcal{L}_{I\beta_I}\{W, Z, S, Y, \beta_0, \theta_0(Z)\} + \mathcal{L}_{I\theta}\{W, Z, S, Y, \beta_0, \theta_0(Z)\}\theta_{\beta_I}(Z, \beta_0)],$$
$$B = \mathrm{cov}[\mathcal{L}_I\{W, Z, S, Y, \beta_0, \theta_0(Z)\} - \mathcal{U}_I(Z)\Psi\{W, Z, S, Y, \beta_0, \theta_0(Z)\}],$$

we obtain the following results.

**Theorem 4.** *Suppose that condition (*P1*) holds. Under regularity conditions (*B1*)–(*B4*), if $\lambda_n \to 0$, $d_n^5 n^{-1} \to 0$, and (2.6) holds, then with probability tending to one, any root $n$ consistent estimator $\hat{\beta}_n = (\hat{\beta}_I^{\mathrm{T}}, \hat{\beta}_{II}^{\mathrm{T}})^{\mathrm{T}}$ obtained in (3.9) must satisfy that*

(i) $\hat{\beta}_{II} = 0$,
(ii) *for any $d_1 \times 1$ vector $v$ such that $v^{\mathrm{T}} v = 1$,*

$$\sqrt{n/d_n} v^{\mathrm{T}} B^{-1/2} (A - \Sigma) \{\hat{\beta}_I - \beta_{I0} - (A - \Sigma)^{-1} b\} \xrightarrow{\mathrm{D}} N(0, 1).$$

The proof of Theorem 4 is given in the Appendix. Theorem 4 implies that for fixed and finite $d$, the convergence rate of the resulting estimate is $n^{-1/2} + a_n$. It also implies that any root $n$ consistent solution $\hat{\beta} = (\hat{\beta}_I^{\mathrm{T}}, \hat{\beta}_{II}^{\mathrm{T}})^{\mathrm{T}}$ of (3.9) must satisfy $\hat{\beta}_{II} = 0$, and $\hat{\beta}_I$ has the following asymptotic normality:

$$\sqrt{n} \{\hat{\beta}_I - \beta_{I0} - (A - \Sigma)^{-1} b\} \xrightarrow{\mathrm{D}} N\{0, (A - \Sigma)^{-1} B (A - \Sigma)^{-\mathrm{T}}\}.$$

See the earlier version of this work, Ma and Li (2007), for details.

## 4. Numerical studies and application

In this section, we provide implementation details such as tuning parameter selection and model error approximation. Issues related to the numerical procedure to solve (2.2) and (3.9), the choice of kernel and bandwidth in the semi-parametric model, and the treatment of multiple roots have been addressed in Ma and Carroll (2006) and are not further discussed here. We assess the finite sample performance of the proposed procedure by Monte Carlo simulation and illustrate the proposed methodology by an empirical analysis of the Framingham heart study data. In our simulation, we concentrate on the performance of the proposed procedure for a quadratic logistic measurement error model and a partially linear logistic measurement error model in terms of model complexity and model error.

### 4.1. Tuning parameter selection

An MM algorithm (Hunter and Li (2005)) and a local linear approximation (LLA) algorithm (Zou and Li (2008)) have been proposed for penalized likelihood with non-concave penalty. However, both the minorize–maximize (MM) algorithm and the LLA algorithm are difficult to implement for the measurement error models we consider. Thus, we employ the local quadratic approximation (LQA) algorithm (Fan and Li (2001)) to solve the penalized estimating equations. Specifically, in implementing the Newton–Raphson algorithm to solve the penalized estimating equations, we locally approximate the first-order derivative of the penalty function by a linear function, following the idea of the LQA. Specifically, suppose that at the $k$th step of the iteration, we obtain the value $\beta^{(k)}$. Then, for $\beta_j^{(k)}$ not very close to zero,

$$p'_\lambda(\beta_j) = p'_\lambda(|\beta_j|)\operatorname{sign}(\beta_j) \approx \frac{p'_\lambda(|\beta_j^{(k)}|)}{|\beta_j^{(k)}|}\beta_j.$$

Otherwise, we set $\beta_j^{(k+1)} = 0$, and exclude the corresponding covariate from the model. This approximation is updated in every step of the Newton–Raphson algorithm iteration. In practice, we set the initial value of $\beta$ to be the unpenalized estimating equation estimate. It can be shown that when the algorithm converges, the solution will satisfy the penalized estimating equations. Following Theorems 2 and 4, we can further approximate the estimation variance of the resulting estimator. That is

$$\widehat{\operatorname{cov}}(\hat{\beta}) = \frac{1}{n}(E - \Sigma_\lambda)^{-1} F (E - \Sigma_\lambda)^{-\mathrm{T}},$$

where $\Sigma_\lambda$ is a diagonal matrix with elements equal to $p'_\lambda(|\hat{\beta}_j|)/|\hat{\beta}_j|$ for non-vanishing $\hat{\beta}_j$, a linear approximation of $\Sigma$ defined in (2.5). We use $E$ to denote the sample approximation of $E\,\partial S^*_{e\!f\!f,I}(W,Z,Y,\beta_I)/\partial\beta_I$ evaluated at $\hat{\beta}$ for the parametric model (2.1) and the sample approximation of the matrix $A$ evaluated at $\hat{\beta}$ for the semi-parametric model (3.8). Similarly, we use $F$ to denote the sample approximation of $\operatorname{cov}(S^*_{e\!f\!f,I})$ evaluated at $\hat{\beta}$ for the parametric model and the sample approximation of the matrix $B$ evaluated at $\hat{\beta}$ for the semi-parametric model, respectively. The consistency of the proposed sandwich formula can be shown by using similar techniques as in Fan and Peng (2004). The accuracy of this sandwich formula will be tested in our simulation studies.

It is desirable to have automatic, data-driven methods to select the tuning parameter $\lambda$. Here we will consider two tuning parameter selectors, the GCV and BIC. To define the GCV and BIC statistics, we need to define the degrees of freedom and goodness-of-fit measure for the final selected model. Similar to the nonconcave penalized likelihood approach, we may define the effective number of parameters or degrees of freedom to be

$$df_\lambda = \operatorname{trace}\{I(I + \Sigma_\lambda)^{-1}\},$$

where $I$ stands for the Fisher information matrix. For the logistic regression models employed in this section, a natural approximation of $I$, ignoring the measurement error effect, is $V^\mathrm{T} Q V$, where $V$ represents the covariates included in the model and $Q$ is a diagonal matrix with the $i$th element equal to $\hat{\mu}_{\lambda,i}(1 - \hat{\mu}_{\lambda,i})$. Here, $\hat{\mu}_{\lambda,i} = P(Y_i = 1|V_i)$.

In the logistic regression model context of this section, we may employ its deviance as a goodness-of-fit measure. Specifically, let $\mu_i$ be the conditional expectation of $Y_i$ given its covariates for $i = 1, \ldots, n$. The deviance of a model fit $\hat{\mu}_\lambda = (\hat{\mu}_{\lambda,1}, \ldots, \hat{\mu}_{\lambda,n})^\mathrm{T}$ is defined to be

$$D(\hat{\mu}_\lambda) = 2\sum_{i=1}^n [Y_i \log(Y_i/\hat{\mu}_{\lambda,i}) + (1 - Y_i)\log\{(1 - Y_i)/(1 - \hat{\mu}_{\lambda,i})\}].$$

Define the GCV statistic to be

$$GCV(\lambda) = \frac{D(\hat{\mu}_\lambda)}{n(1 - df_\lambda/n)^2},$$

and the BIC statistic to be

$$BIC(\lambda) = D(\hat{\mu}_\lambda) + 2\log(n)df_\lambda.$$

The GCV and the BIC tuning parameter selectors select $\lambda$ by minimizing $GCV(\lambda)$ and $BIC(\lambda)$, respectively. Note that the BIC tuning parameter selector is distinguished from the traditional BIC variable selection criterion, which is not well defined for estimating equation methods. Wang, Li and Tsai (2007) provided a study on the asymptotic behavior for the GCV and BIC tuning parameter selectors for the non-concave penalized least-squares variable selection procedures in linear and partially linear regression models. Further study of the asymptotic property of the proposed tuning parameter selection is needed, but it is outside the scope of this paper.

### 4.2. Model error

Model error is an effective way of evaluating model adequacy versus model complexity. To implement the concept of model error in evaluating our procedure, we first simplify its definition for the logistic partially linear measurement error model. Denote $\mu(S, X, Z) = E(Y|S, X, Z)$, and define the model error for a model $\hat{\mu}(S, X, Z)$ as

$$ME(\hat{\mu}) = E\{\hat{\mu}(S^+, X^+, Z^+) - \mu(S^+, X^+, Z^+)\}^2,$$

where the expectation is taken over the new observation $S^+$, $X^+$ and $Z^+$. Let $g(\cdot)$ be the logit link. For the logistic partially linear model, the mean function has the form $\mu(S, X, Z) = g^{-1}\{\theta(Z) + \beta^T V\}$, where $V = (S^T, X^T)^T$. If $\hat{\theta}(\cdot)$ and $\hat{\beta}$ are consistent estimators for $\theta(\cdot)$ and $\beta$, respectively, then by a Taylor expansion the model error can be approximated by

$$ME(\hat{\mu}) \approx E(\dot{g}^{-1}\{\theta(Z^+) + \beta^T V^+\}^2[\{\hat{\theta}(Z^+) - \theta(Z^+)\}^2 \\ + (\hat{\beta}^T V^+ - \beta^T V^+)^2 + 2\{\hat{\theta}(Z^+) - \theta(Z^+)\}(\hat{\beta}^T V^+ - \beta^T V^+)]).$$

The first component is the inherent model error due to $\hat{\theta}(\cdot)$, the second one is due to the lack of fit of $\hat{\beta}$, and the third one is the cross-product between the first two components. Thus, to assess the performance of the proposed variable selection procedure, we define the approximate model error (AME) for $\hat{\beta}$ to be

$$AME(\hat{\beta}) = E[\dot{g}^{-1}\{\theta(Z^+) + \beta^T V^+\}^2(\hat{\beta}^T V^+ - \beta^T V^+)^2].$$

Furthermore, the AME of $\hat{\beta}$ can be written as

$$AME(\hat{\beta}) = (\hat{\beta} - \beta)^{\mathrm{T}} E[\dot{g}^{-1}\{\theta(Z^+) + \beta^{\mathrm{T}} V^+\}^2 V^+ V^{+\mathrm{T}}](\hat{\beta} - \beta) \qquad (4.1)$$
$$\hat{=} (\hat{\beta} - \beta)^{\mathrm{T}} C_X (\hat{\beta} - \beta).$$

In our simulation, the matrix $C_X$ is estimated by 1 million Monte Carlo simulations. For measurement error data, we observe $W$ instead of $X$. We also consider an alternative approximate model error

$$AME_W(\hat{\beta}) = (\hat{\beta} - \beta)^{\mathrm{T}} C_W (\hat{\beta} - \beta), \qquad (4.2)$$

where $C_W$ is obtained by replacing $X$ with $W$ in the definition of $C_X$. The $AME(\hat{\beta})$ and $AME_W(\hat{\beta})$ are defined for the parametric model case by setting $\theta(\cdot) = 0$. Note that although we defined the model error in the context with a logistic link function, it is certainly not restricted to such a case. The general approach for calculating $AME$ is to approximate the probability density function evaluated at the estimated parameters around the true parameter value and to extract the linear term of the parameter of interest. $AME_W$ is calculated by replacing $X$ with $W$.

### 4.3. Simulation examples

To demonstrate the performance of our method in both parametric and semi-parametric measurement error models, we conduct two simulation studies. In our simulation, we will examine only the performance of the penalized estimating equation method with the SCAD penalty.

***Example 1.*** In this example, we generate data from a logistic model where the covariate measured with error enters the model through a quadratic function and the covariates measured without error enter linearly. The measurement error follows a normal additive pattern. Specifically,

$$\mathrm{logit}\{p(Y=1|X,Z)\} = \beta_0 + \beta_1 X + \beta_2 X^2 + \beta_3 Z_1 + \beta_4 Z_2 + \beta_5 Z_3 + \beta_6 Z_4$$
$$+ \beta_7 Z_5 + \beta_8 Z_6 + \beta_9 Z_7$$

and

$$W = X + U,$$

where $\beta = (0, 1.5, 2, 0, 3, 0, 1.5, 0, 0, 0)^{\mathrm{T}}$, the covariate $X$ is generated from a normal distribution with mean 0 and variance 1, $(Z_1, \ldots, Z_6)^{\mathrm{T}}$ is generated from a normal distribution with mean 0, and covariance between $Z_i$ and $Z_j$ is $0.5^{|i-j|}$. The last component of the $Z$ covariates, $Z_7$, is a binary variable taking value 0 or 1 with equal probability. $U$ is normally distributed with mean 0 and standard deviation 0.1. In our simulation, the sample size is taken to be either $n = 500$ or $n = 1000$.

**Table 1.** MRMEs and model complexity, for example, 1

|     | $n$  | RAME median (MAD) | $RAME_W$ median (MAD) | # of zero coefficients C | E |
|-----|------|----------------|----------------|-------|-------|
| GCV | 500  | 0.694 (0.231)  | 0.698 (0.228)  | 4.574 | 0.006 |
| BIC | 500  | 0.396 (0.188)  | 0.396 (0.187)  | 5.857 | 0.074 |
| GCV | 1000 | 0.766 (0.187)  | 0.770 (0.185)  | 4.456 | 0     |
| BIC | 1000 | 0.390 (0.157)  | 0.401 (0.158)  | 5.758 | 0.010 |

For the selected model, the model complexity is summarized in terms of the number of zero coefficients and the model error is summarized in terms of relative approximation model error (RAME), defined to be the ratio of model error of the selected model to that of the full model. In Table 1, the RAME column corresponds to the sample median and median absolute deviation (MAD) divided by a factor of 0.6745 of the RAME values over 1000 simulations. Similarly, the $RAME_W$ column corresponds to those of the $RAME_W$ values over 1000 simulations. From Table 1, it can be seen that the values of RAME and $RAME_W$ are very close. The average count of zero coefficients is also reported in Table 1, where the column labeled "C" presents the average count restricted only to the true zero coefficients, while the column labeled "E" displays the average count of the coefficients erroneously set to 0.

We next verify the consistency of the estimators and test the accuracy of the proposed standard error formula. Table 2 displays the bias and sample standard deviation (SD) of the estimates for two non-zero coefficients, $(\beta_1, \beta_2)$, over 1000 simulations and the sample average and the sample standard deviations of the 1000 standard errors obtained by using the sandwich formula. The row labeled "EE" corresponds to the unpenalized estimating equation estimator. We omit here the results for other non-zero coefficients and the results under sample size $n = 500$. Interested readers can find them in an earlier version of this work, Ma and Li (2007). Overall, the estimators are consistent and the sandwich formula works well.

**Table 2.** Bias and standard errors, for example, 1 ($n = 1000$)

|     | $\hat{\beta}_1$ | | $\hat{\beta}_2$ | |
|-----|--------------|--------------|--------------|--------------|
|     | Bias (SD)    | SE (Std(SE)) | Bias (SD)    | SE (Std(SE)) |
| EE  | 0.072 (0.273) | 0.268 (0.062) | 0.124 (0.332) | 0.321 (0.088) |
| GCV | 0.029 (0.254) | 0.250 (0.048) | 0.009 (0.258) | 0.253 (0.057) |
| BIC | 0.024 (0.290) | 0.249 (0.054) | 0.052 (0.255) | 0.244 (0.052) |

**Table 3.** MRMEs and model complexity, for example, 2

| Method | $n$ | $\text{RME}_X$ | $\text{RME}_W$ | # of zero coefficients | |
|---|---|---|---|---|---|
| | | median (MAD) | median (MAD) | C | E |
| GCV | 500 | 0.878 (0.161) | 0.880 (0.158) | 4.060 | 0 |
| BIC | 500 | 0.381 (0.158) | 0.387 (0.155) | 5.713 | 0 |
| GCV | 1000 | 0.868 (0.164) | 0.873 (0.160) | 4.061 | 0 |
| BIC | 1000 | 0.386 (0.162) | 0.392 (0.161) | 5.694 | 0 |

*Example 2.* In this example, we illustrate the performance of the method for a semiparametric measurement error model. Simulation data are generated from

$$\text{logit}(Y) = \beta_1 X + \beta_2 S_1 + \cdots + \beta_{10} S_9 + \theta(Z),$$
$$W = X + U,$$

where $\beta$, $X$, and $W$ are the same as in the previous simulation. We generate $S$'s in a fashion similar to the $Z$'s in Example 1. That is, $(S_1, \ldots, S_8)$ is generated from a normal distribution with mean zero and covariance between $S_i$ and $S_j$ is $0.5^{|i-j|}$. $S_9$ is a binary variable with equal probability to be zero or one. The random variable $Z$ is generated from a uniform distribution in $[-\pi/2, \pi/2]$. The true function $\theta(Z) = 0.5\cos(Z)$. The parameter takes values $\beta = (1.5, 2, 0, 0, 3, 0, 1.5, 0, 0, 0)^{\text{T}}$.

The simulation results are summarized in Table 3, with notation similar to that of Table 1. From Table 3, we can see that the penalized estimating equation estimators can significantly reduce model complexity. Overall, the BIC tuning parameter selectors perform better, while GCV is too conservative. We have further tested the consistency and the accuracy of the standard error formula derived from the sandwich formula. The result is summarized in Table 4, with notation similar to that of Table 2. We note the consistency of the estimator and that the standard error formula performs very well. More simulation results are summarized in the earlier version of the work, Ma and Li (2007).

**Table 4.** Bias and standard errors, for example, 2 ($n = 1000$)

| Method | $\hat{\beta}_1$ | | $\hat{\beta}_2$ | |
|---|---|---|---|---|
| | Bias (SD) | SE (Std(SE)) | Bias (SD) | SE (Std(SE)) |
| EE | 0.039 (0.170) | 0.166 (0.018) | 0.057 (0.194) | 0.190 (0.018) |
| GCV | 0.047 (0.174) | 0.172 (0.020) | 0.069 (0.196) | 0.191 (0.021) |
| BIC | 0.031 (0.169) | 0.170 (0.019) | 0.044 (0.179) | 0.185 (0.019) |

### 4.4. An application

The Framingham heart study data set (Kannel *et al.* (1986)) is a well-known data set where it is generally accepted that measurement error exists on the long-term systolic blood pressure (SBP). In addition to SBP, other measurements include age, smoking status, and serum cholesterol. In the literature, there has been speculation that a second-order term involving age might be needed to analyze the dependence of heart disease occurrence. In addition, it is unclear if the interaction between the various covariates plays a role in influencing the heart disease rate. The data set includes 1615 observations.

With the method developed here, it is possible to perform a variable selection to address these issues. Following the literature, we adopt the measurement error model of $\log(\text{MSBP} - 50) = \log(\text{SBP} - 50) + U$, where $U$ is a mean zero normal random variable with variance $\sigma_u^2 = 0.0126$ and MSBP is the measured SBP. We denote the standardized $\log(\text{MSBP} - 50)$ as $W$. The standardization using the same parameters on $\log(\text{SBP} - 50)$ is denoted $X$. The standardized serum cholesterol and age are denoted by $Z_1, Z_2$, respectively, and $Z_3$ denotes the binary variable smoking status. Using $Y$ to denote the occurrence of heart disease, the saturated model that includes all the interaction terms and also the square of age term is of the form

$$\text{logit}\{p(Y = 1|X, Z's)\} = \beta_1 X + \beta_2 X Z_1 + \beta_3 X Z_2 + \beta_4 X Z_3 + \beta_5 + \beta_6 Z_1 + \beta_7 Z_2$$
$$+ \beta_8 Z_3 + \beta_9 Z_2^2 + \beta_{10} Z_1 Z_2 + \beta_{11} Z_1 Z_3 + \beta_{12} Z_2 Z_3,$$
$$W = X + U.$$

We used both GCV and BIC tuning parameter selectors to choose $\lambda$. We present the tuning parameters and the corresponding GCV and BIC scores in Figure 1. The final chosen $\lambda$ is 0.073 and 0.172 by the GCV and BIC selectors, respectively. The selected model is depicted in Table 5. The GCV criterion selects the covariates $X, XZ_1, 1, Z_1, Z_2, Z_3, Z_2^2, Z_2 Z_3$ into the model, while the BIC criterion selects the covariates $X, 1, Z_1, Z_2$ into the model. We report the selection and estimation results in Table 5, as well as the semi-parametric estimation results without variable selection.

As shown, the terms $X, 1, Z_1, Z_2$ are selected by both criteria, while $Z_3, Z_2^2$, and some of the interaction terms are selected only by GCV. The BIC criterion is very aggressive and it results in a very simple final model while the GCV criterion is much more conservative, hence the resulting model is more complex. This agrees with the simulation results obtained. Since both criteria have included the covariate $X$, the measurement error feature and its treatment in the Framingham data is unavoidable.

## 5. Discussion

In this paper, we have proposed a new class of variable selection procedures in the framework of measurement error models. The procedure is proposed in a completely general functional measurement error model setting and is suitable for both parametric and semi-parametric models that contain unspecified smooth functions of an observable covariate.

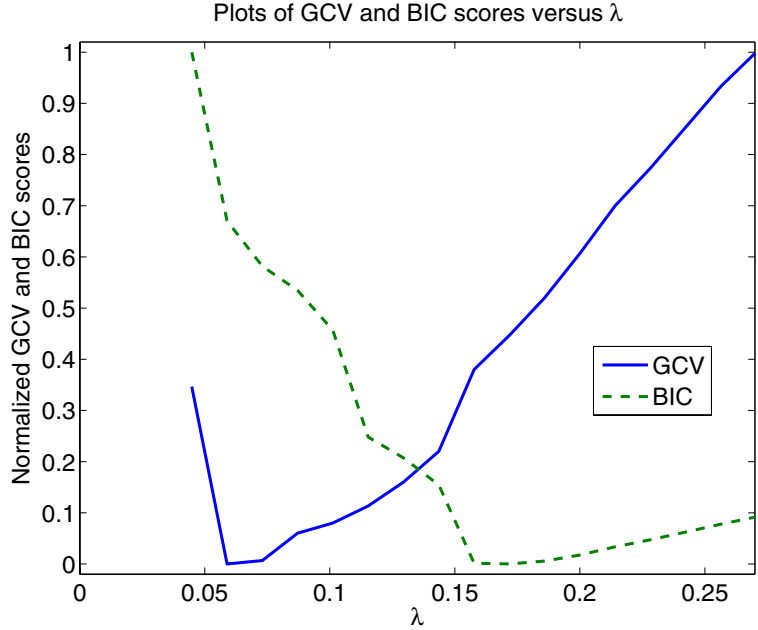

**Figure 1.** Tuning parameters and their corresponding BIC and GCV scores for the Framingham data. The scores are normalized to the range $[0, 1]$.

We have assumed the error model $p_{W|X,Z}(W|X, Z)$ to be completely known for ease of presentation. When the error model contains an unknown parameter $\xi$, the identi-fiability of the problem requires additional information such as multiple measurements or instruments. Such information should be incorporated to estimate $\xi$. Specifically, in the variable selection context, we can simply append the estimating equation with these additional estimating equations obtained from the corresponding score functions with respect to $\xi$ and append the penalty function $p'_\lambda$ with zeros. Because the augmented estimating equations still have the same convergence property as illustrated in Ma and Carroll (2006), the same asymptotic convergence rates and oracle properties hold as in the known $\xi$ case, without any efficiency loss. When the error model $p_{W|X,Z}(W|X,Z)$ is completely unspecified, a nonparametric estimation of the measurement error distribution has to be carried out first, then the result can be plugged into the proposed variable selection and estimation procedure. In this case, the asymptotic convergence rate of the parameters and the oracle property remain the same, but the asymptotic variance will increase. The details of incorporating the estimation of unknown error and demonstrating its subsequent convergence property in the estimation framework are the focus of Hall and Ma (2007).

**Table 5.** Results for the Framingham data set

|  | EE | GCV | BIC |
|---|---|---|---|
|  | $\hat{\beta}$ (SE) | $\hat{\beta}$ (SE) | $\hat{\beta}$ (SE) |
| $X$        | 0.643 (0.248)  | 0.416 (0.093)  | 0.179 (0.039) |
| $XZ_1$     | $-0.167$ (0.097) | $-0.072$ (0.041) | 0 (NA) |
| $XZ_2$     | $-0.059$ (0.111) | 0 (NA)         | 0 (NA) |
| $XZ_3$     | $-0.214$ (0.249) | 0 (NA)         | 0 (NA) |
| Intercept  | $-3.415$ (0.428) | $-3.255$ (0.356) | $-2.555$ (0.092) |
| $Z_1$      | 0.516 (0.212)  | 0.332 (0.085)  | 0.124 (0.033) |
| $Z_2$      | 1.048 (0.341)  | 1.044 (0.329)  | 0.398 (0.067) |
| $Z_3$      | 1.060 (0.443)  | 0.907 (0.373)  | 0 (NA) |
| $Z_2^2$    | $-0.253$ (0.125) | $-0.262$ (0.121) | 0 (NA) |
| $Z_1 Z_2$  | $-0.072$ (0.103) | 0 (NA)         | 0 (NA) |
| $Z_1 Z_3$  | $-0.161$ (0.225) | 0 (NA)         | 0 (NA) |
| $Z_2 Z_3$  | $-0.442$ (0.336) | $-0.473$ (0.326) | 0 (NA) |

We also would like to point out that in the special case of generalized linear models and normal additive error with possible heteroscedasticity, the procedure of solving linear integral equations can be spared and the estimating equations are simplified significantly (Ma and Tsiatis (2006)). In such situations, the computation complexity of the proposed procedure will be reduced to about the same level as for variable selection in regressions without errors in the variables.

As pointed out by the referee, it is interesting to perform variable selection for high-dimensional data. In this paper, we allow the number of covariates to grow to infinity at a $o_p(n^{-1/5})$ rate as the sample size $n$ increases. However, the proposed procedures and the used algorithm in this paper may not be directly applied to the large $p$, small $n$ problems. Variable selection for the large $p$, small $n$ setting is a very active research topic. It is challenging to extend the existing variable selection procedures for large $p$, small $n$ problems to measurement error data. Further research is needed on this topic, but this is outside the scope of this paper.

## Appendix

*Global assumption* (P1) *on the penalty function:*

(P1) Let $c_n = \max\{|p''_\lambda(|\beta_{j0}|)|: \beta_{j0} \neq 0\}$. Assume that $\lambda_n \to 0$, $a_n = \mathrm{O}(n^{-1/2})$ and $c_n \to 0$ $n \to \infty$. In addition, there exist constants $C$ and $D$ such that when $\gamma_1, \gamma_2 > C\lambda$, $|p''_\lambda(\gamma_1) - p''_\lambda(\gamma_2)| \leq D|\gamma_1 - \gamma_2|$.

It is easy to verify that both the $L_1$ and the SCAD penalties satisfy this condition.

*The regularity conditions for Theorems* 1 *and* 2:

(A1) The expectation of the first derivative of $S^*_{\textit{eff}}$ with respect to $\beta$ exists at $\beta_0$ and its left eigenvalues are bounded away from zero and infinity uniformly for all $n$. For any entry $S_{jk}$ in $\partial S^*_{\textit{eff}}(\beta_0)/\partial \beta^{\mathrm{T}}$, $E(S^2_{jk}) < C_1 < \infty$.

(A2) The eigenvalues of the matrix $E(S^*_{\textit{eff},I} S^{*\mathrm{T}}_{\textit{eff},I})$ satisfy $0 < C_2 < \lambda_{\min} < \cdots < \lambda_{\max} < C_3 < \infty$ for all $n$. For any entries, $S_k, S_j$ in $S^*_{\textit{eff}}(\beta_0)$, $E(S_k^2 S_j^2) < C_4 < \infty$.

(A3) The second derivatives of $S^*_{\textit{eff}}$ with respect to $\beta$ exist and the entries are uniformly bounded by a function $M(W_i, Z_i, Y_i)$ in a large enough neighborhood of $\beta_0$. In addition, $E(M^2) < C_5 < \infty$ for all $n, d$.

Conditions (A1)–(A3) are mild regularity conditions. They guarantee that the solution of the following estimating equation

$$\sum_{i=1}^n S^*_{\textit{eff}}(W_i, Z_i, Y_i, \beta) = 0$$

is root $(n/d_n)$ convergent, and possesses asymptotic normality.

**Proof of Theorem 1.** Condition (A1) allows us to define

$$J = \left\{ E\left(\left.\frac{\partial S^*_{\textit{eff}}}{\partial \beta^{\mathrm{T}}}\right|_{\beta_0}\right)\right\}^{-1}, \qquad \phi^*_{\textit{eff}} = J S^*_{\textit{eff}} \quad \text{and} \quad q'_{\lambda_n}(\beta) = J p'_{\lambda_n}(\beta).$$

Let $\alpha_n = n^{-1/2} + a_n$ and $\phi^*_{\textit{eff},i}(\beta) = \phi^*_{\textit{eff}}(W_i, Z_i, Y_i, \beta)$. It suffices to show that

$$n^{-1/2} \sum_{i=1}^n \phi^*_{\textit{eff},i}(\beta) - n^{1/2} q'_{\lambda_n}(\beta) = 0 \tag{A1}$$

has a solution $\hat{\beta}$ that satisfies $\|\hat{\beta} - \beta_0\| = \mathrm{O}_p(\sqrt{d_n}\alpha_n)$. This will be shown using the Brouwer fixed point theorem. Using the Taylor expansion, we have

$$n^{-1/2} \sum_{i=1}^n \phi^*_{\textit{eff},i}(\beta) - n^{1/2} q'_{\lambda_n}(\beta)$$

$$= n^{-1/2} \sum_{i=1}^n \phi^*_{\textit{eff},i}(\beta_0) - n^{1/2} q'_{\lambda_n}(\beta_0)$$

$$+ n^{-1/2} \sum_{i=1}^n \frac{\partial \phi^*_{\textit{eff},i}(\beta^*)}{\partial \beta^{\mathrm{T}}}(\beta - \beta_0)$$

$$- n^{1/2} \frac{\partial q'_{\lambda_n}(\beta_0)}{\partial \beta^{\mathrm{T}}}(\beta - \beta_0)\{1 + \mathrm{o}_p(1)\},$$

where $\beta^*$ is between $\beta$ and $\beta_0$. It can be shown by conditions (A1)–(A3) and definition of $\phi^*_{\textit{eff}}(\cdot)$ that

$$(\beta - \beta_0)^{\mathrm{T}} \left\{ \frac{1}{n} \sum_{i=1}^n \frac{\partial \phi^*_{\textit{eff},i}(\beta^*)}{\partial \beta^{\mathrm{T}}} \right\}(\beta - \beta_0) = \|\beta - \beta_0\|^2 \{1 + \mathrm{o}_P(1)\}.$$

We next check the key condition for the Brouwer fixed point theorem. For any $\beta$ such that $\|\beta - \beta_0\| = C\sqrt{d_n}\alpha_n$ for some constant $C$, it follows by condition (P1) that

$$(\beta - \beta_0)^{\mathrm{T}}\left\{\frac{1}{\sqrt{n}}\sum_{i=1}^{n}\phi^*_{\mathit{eff},i}(\beta) - n^{1/2}q'_{\lambda_n}(\beta)\right\}$$

$$= (\beta - \beta_0)^{\mathrm{T}}\left\{\frac{1}{\sqrt{n}}\sum_{i=1}^{n}\phi^*_{\mathit{eff},i}(\beta_0) - \sqrt{n}q'_{\lambda_n}(\beta_0)\right\} + \sqrt{n}\|\beta - \beta_0\|^2\{1 + \mathrm{o}_P(1)\}.$$

Using the Cauchy–Schwarz inequality, it can be shown that the first term in the above equation is of order $\|\beta - \beta_0\|\mathrm{O}_p(\sqrt{d_n + d_n n \alpha_n^2}) = \mathrm{O}_p(Cn^{1/2}d_n\alpha_n^2)$. Note that $\sqrt{n}\|\beta - \beta_0\|^2 = C^2 n^{1/2}d_n\alpha_n^2$. Thus the second term in the above equation dominates the first term with probability $1 - \epsilon$ for any $\epsilon > 0$ as long as $C$ is large enough. Thus, for any $\epsilon > 0$ and for large enough $C$, the probability for the above display to be larger than zero is at least $1 - \epsilon$. From the Brouwer fixed point theorem, we know that with a probability of at least $1 - \epsilon$, there exists at least one solution for (A1) in the region $\|\beta - \beta_0\| \leq C\sqrt{d_n}\alpha_n$. $\square$

*Lemma on sparsity for Theorem* 2.

**Lemma 1.** *If the conditions in Theorem 2 hold, then for any given $\beta$ that satisfies $\|\beta - \beta_0\| = \mathrm{O}_p(\sqrt{d_n/n})$, with probability tending to 1, any solution $(\beta_I^{\mathrm{T}}, \beta_{II}^{\mathrm{T}})^{\mathrm{T}}$ of (2.2) satisfies that $\beta_{II} = 0$.*

**Proof.** Denote the $k$th element in $\sum_{i=1}^{n}S^*_{\mathit{eff}}(W_i, Z_i, Y_i, \beta)$ as $L_{nk}(\beta)$, $k = d_1 + 1, \ldots, d_n$. We next show that the order of $L_{nk}(\beta)$ is $\mathrm{O}_p(\sqrt{nd_n})$,

$$L_{nk}(\beta) = L_{nk}(\beta_0) + \sum_{j=1}^{d_n}\frac{\partial L_{nk}(\beta_0)}{\partial \beta_j}(\beta_j - \beta_{j0}) \qquad (A2)$$

$$+ 2^{-1}\sum_{l=1}^{d_n}\sum_{j=1}^{d_n}\frac{\partial^2 L_{nk}(\beta^*)}{\partial \beta_l \partial \beta_j}(\beta_l - \beta_{l0})(\beta_j - \beta_{j0}),$$

where $\beta^*$ is between $\beta$ and $\beta_0$. Because of condition (A2), the first term of (A2) is of order $\mathrm{O}_p(n^{1/2}) = \mathrm{o}_p(\sqrt{nd_n})$. The second term in (A2) can be further written as

$$\sum_{j=1}^{d_n}\left\{\frac{\partial L_{nk}(\beta_0)}{\partial \beta_j} - E\frac{\partial L_{nk}(\beta_0)}{\partial \beta_j}\right\}(\beta_j - \beta_{j0}) + \sum_{j=1}^{d_n}E\frac{\partial L_{nk}(\beta_0)}{\partial \beta_j}(\beta_j - \beta_{j0}). \qquad (A3)$$

Using the Cauchy–Schwarz inequality and condition (A1), it can be shown by straightforward calculation that the first term in (A3) is of order $\mathrm{O}_p(\sqrt{d_n/n}) = \mathrm{o}_p(\sqrt{nd_n})$. Using the Cauchy–Schwarz inequality again, the second term in (A3) is controlled by

$$n\left\{\sum_{j=1}^{d_n}\left(E\frac{\partial S_{\textit{eff},k}}{\partial \beta^{\mathrm{T}}}\right)^2\right\}^{1/2}\|\beta-\beta_0\| \leq n\lambda_{\max}\left(E\frac{\partial S_{\textit{eff},k}}{\partial \beta^{\mathrm{T}}}\right)^2\|\beta-\beta_0\| = \mathrm{O}_p(\sqrt{nd_n}).$$

Thus, the second term of (A2) is of order $\mathrm{O}_p(\sqrt{nd_n})$. As for the third term of (A2), we can have a similar decomposition to that of (A3). Using the Cauchy–Schwarz inequality in matrix form and condition (A3), it can be shown that the third term of (A2) is of order $\mathrm{O}_p(d_n^2) + \mathrm{O}_p(n^{-1/2}d_n^2) = \mathrm{o}_p(\sqrt{nd_n})$ as $d_n^5/n \to 0$. Thus, $L_{nk}(\beta)$ is of order $\mathrm{O}_p(\sqrt{nd_n})$. Hence we have

$$L_{nk}(\beta) - np'_{\lambda_n}(\beta_k) = -\sqrt{n}\{\sqrt{n/d_n}p'_{\lambda_n}(|\beta_k|)\,\mathrm{sign}(\beta_k) + \mathrm{O}_p(1)\}.$$

Using condition (2.6), the sign of $L_{nk}(\beta) - np'_{\lambda_n}(\beta_k)$ is decided by $\mathrm{sign}(\beta_k)$ completely when $n$ is large enough. From the continuity of $L_{nk}(\beta) - np'_{\lambda_n}(\beta_k)$, we obtain that it is zero at $\beta_k = 0$. □

**Proof of Theorem 2.** From Theorem 1 and condition (P1), there is a root $(n/d_n)$ consistent estimator $\hat{\beta}$. From Lemma 1, $\hat{\beta} = (\hat{\beta}_I^{\mathrm{T}}, 0^{\mathrm{T}})^{\mathrm{T}}$, so (i) is shown. Denote the first $d_1$ equations in $\sum_{i=1}^n S_{\textit{eff}}^*\{W_i, Z_i, Y_i, (\beta_I^{\mathrm{T}}, 0^{\mathrm{T}})^{\mathrm{T}}\}$ as $L_n(\beta_I)$. Now consider solving the first $d_1$ equations in (2.2) for $\beta_I$, while $\beta_{II} = 0$. We have

$$0 = L_n(\hat{\beta}_I) - np'_{\lambda_n,I}(\hat{\beta}_I)$$
$$= L_n(\beta_{I0}) + \frac{\partial L_n(\beta_{I0}^*)}{\partial \beta_I^{\mathrm{T}}}(\hat{\beta}_I - \beta_{I0}) - nb_n - np''_{\lambda_n}(\beta_I^*)(\hat{\beta}_I - \beta_{I0}),$$

where $\beta_I^*$ is between $\beta_{I0}$ and $\hat{\beta}_I$. It follows by condition (P1) that

$$\left\|n^{-1}\frac{\partial L_n(\beta_{I0}^*)}{\partial \beta_I^{\mathrm{T}}} - p''_{\lambda_n,I}(\beta_I^*) - E\frac{\partial L_n(\beta_{I0})}{\partial \beta_I^{\mathrm{T}}} + p''_{\lambda_n,I}(\beta_{I0})\right\|^2$$
$$\leq 2\left\|n^{-1}\frac{\partial L_n(\beta_{I0}^*)}{\partial \beta_I^{\mathrm{T}}} - E\frac{\partial L_n(\beta_{I0})}{\partial \beta_I^{\mathrm{T}}}\right\|^2 + \mathrm{O}_p(n^{-1}d_n).$$

Furthermore, for any fixed $\epsilon > 0$, it follows by conditions (A1) and (A3) and the Chebyshev inequality that

$$P_r\left\{\left\|n^{-1}\frac{\partial L_n(\beta_{I0}^*)}{\partial \beta_I^{\mathrm{T}}} - E\frac{\partial L_n(\beta_{I0})}{\partial \beta_I^{\mathrm{T}}}\right\| \geq \epsilon d_n^{-1}\right\}$$
$$\leq \frac{d_n^2}{n^2\epsilon^2}E\left\|\frac{\partial L_n(\beta_{I0}^*)}{\partial \beta_I^{\mathrm{T}}} - nE\frac{\partial L_n(\beta_{I0})}{\partial \beta_I^{\mathrm{T}}}\right\|^2 = \mathrm{O}(d_n^2 n^{-2}d_1^2 n) = \mathrm{o}(1),$$

since $d_1 \leq d_n$. Thus, we have

$$\left\| n^{-1} \frac{\partial L_n(\beta_{I0}^*)}{\partial \beta_I^T} - E \frac{\partial L_n(\beta_{I0})}{\partial \beta_I^T} \right\| = \mathrm{o}_p(d_n^{-1}).$$

Therefore,

$$\left\| n^{-1} \frac{\partial L_n(\beta_{I0}^*)}{\partial \beta_I^T} - p''_{\lambda_n, I}(\beta^*) - E \frac{\partial L_n(\beta_{I0})}{\partial \beta_I^T} + p''_{\lambda_n, I}(\beta_{I0}) \right\|^2 = \mathrm{o}_p(d_n^{-2}),$$

and subsequently,

$$\left\| \left\{ n^{-1} \frac{\partial L_n(\beta_{I0}^*)}{\partial \beta_I^T} - p''_{\lambda_n, I}(\beta^*) - E \frac{\partial L_n(\beta_{I0})}{\partial \beta_I^T} + p''_{\lambda_n, I}(\beta_{I0}) \right\}(\hat{\beta}_I - \beta_{I0}) \right\|$$
$$\leq \mathrm{o}_p(d_n^{-1}) \mathrm{O}_p(n^{-1/2} d_n^{1/2}) = \mathrm{o}_p(n^{-1/2}).$$

We thus obtain that

$$\left\{ -E \frac{\partial L_n(\beta_{I0})}{\partial \beta_I^T} + \Sigma_n \right\}(\hat{\beta}_I - \beta_{I0}) + b_n = n^{-1} L_n(\beta_{I0}) + \mathrm{o}_p(n^{-1/2}).$$

Denote $I^* = E\{S^*_{n,\mathit{eff},I}(\beta_{I0}) S^{*T}_{n,\mathit{eff},I}(\beta_{I0})\}$. Using condition (A2), it follows that

$$n^{1/2} v^T I^{*-1/2} \left[ \left\{ -E \frac{\partial L_n(\beta_{I0})}{\partial \beta_I^T} + \Sigma_n \right\}(\hat{\beta}_I - \beta_{I0}) + b_n \right] = n^{-1/2} v^T I^{*-1/2} L_n(\beta_{I0}) + \mathrm{o}_p(1).$$

Let $Y_i = n^{-1/2} v^T I^{*-1/2} S_{n,\mathit{eff},I}(W_i, Z_i, Y_i, \beta_{I0})$. It follows that for any $\epsilon > 0$,

$$\sum_{i=1}^n E\|Y_i\|^2 \mathbf{1}(\|Y_i\| > \epsilon) = nE\|Y_1\|^2 \mathbf{1}(\|Y_1\| > \epsilon) \leq n(E\|Y_1\|^4)^{1/2} \{P_r(\|Y_1\| > \epsilon)\}^{1/2}.$$

Using Chebyshev's inequality, we have

$$P_r(\|Y_1\| > \epsilon) \leq \frac{E\|Y_1\|^2}{\epsilon^2} = \frac{E\|vI^{*-1/2} S_{\mathit{eff},I}(W_1, Z_1, Y_1, \beta_{I0})\|^2}{n\epsilon^2} = \frac{v^T v}{n\epsilon^2} = \mathrm{O}(n^{-1}).$$

Note that $E(\|Y_1\|^4) = n^{-2} E\|v^T I^{*-1/2} S_{\mathit{eff},I}(W_1, Z_1, Y_1, \beta_{I0})\|^4$. Note that the rank of $vv^T$ is one, and hence $\lambda_{\max}(vv^T)$ equals the trace of $vv^T$. So, $\lambda_{\max}(vv^T) = 1$ as $v^T v = 1$. Thus, it follows that

$$E(\|Y_1\|^4) = n^{-2} E\{S_{\mathit{eff},I}(W_1, Z_1, Y_1, \beta_{I0})^T I^{*-1/2} vv^T I^{*-1/2} S_{\mathit{eff},I}(W_1, Z_1, Y_1, \beta_{I0})\}^2$$
$$\leq n^{-2} \lambda_{\max}^2(I^{*-1}) E\{S_{\mathit{eff},I}(W_1, Z_1, Y_1, \beta_{I0})^T S_{\mathit{eff},I}(W_1, Z_1, Y_1, \beta_{I0})\}^2$$
$$= n^{-2} \lambda_{\max}^2(I^{*-1}) E\|S_{\mathit{eff},I}(W_1, Z_1, Y_1, \beta_{I0})\|^4 = \mathrm{O}(d_1^2 n^{-2}),$$

due to condition (A2). Hence,

$$\sum_{i=1}^{n} E\|Y_i\|^2 \mathbf{1}(\|Y_i\| > \epsilon) = \mathrm{O}(nd_1 n^{-1} n^{-1/2}) = \mathrm{o}(1).$$

On the other hand,

$$\sum_{i=1}^{n} \mathrm{cov}(Y_i) = n\, \mathrm{cov}\{n^{-1/2} v I^{*-1/2} S_{\mathit{eff},I}(W_1, Z_1, Y_1, \beta_{I0})\}$$

$$= v I^{*-1/2} E\{S_{\mathit{eff},I}(W_1, Z_1, Y_1, \beta_{I0}) S_{\mathit{eff},I}(W_1, Z_1, Y_1, \beta_{I0})^{\mathrm{T}}\} I^{*-1/2} v^{\mathrm{T}} = 1.$$

Following the Lindeberg–Feller central limit theorem, the results in (ii) now follow. □

*Regularity conditions for Theorems 3 and 4.*

The notation $C_i$ below is generic and is allowed to be different from that in conditions (A1)–(A3).

(B1) The first derivatives of $\mathcal{L}$ with respect to $\beta$ and $\theta$ exist and are denoted as $\mathcal{L}_\beta$ and $\mathcal{L}_\theta$, respectively. The first derivative of $\theta$ with respect to $\beta$ exists and is denoted as $\theta_\beta$. Thus $E(\mathcal{L}_\beta + \mathcal{L}_\theta \theta_\beta)$ exists and its left eigenvalues are bounded away from zero and infinity uniformly for all $n$ at $\beta_0$ and the true function $\theta_0(Z)$. For any entry $S_{jk}$ of the matrix $\mathrm{d}(\mathcal{L}_\beta + \mathcal{L}_\theta \theta_\beta)/\mathrm{d}\beta$, $E(S_{jk}^2) < C_1 < \infty$.

(B2) The eigenvalues of the matrix $E\{\mathcal{L}_I - \mathcal{U}_I(Z)\Psi\}\{\mathcal{L}_I - \mathcal{U}_I(Z)\Psi\}^{\mathrm{T}}$ satisfy $0 < C_2 < \lambda_{\min} < \cdots < \lambda_{\max} < C_3 < \infty$ for all $n$; for any entries $S_k, S_j$ in $(\mathcal{L}_\beta + \mathcal{L}_\theta \theta_\beta)$, $E(S_k^2 S_j^2) < C_4 < \infty$.

(B3) The second derivatives of $\mathcal{L}$ with respect to $\beta$ and $\theta$ exist, the second derivatives of $\theta$ with respect to $\beta$ exist, and the entries are uniformly bounded by a function $M(W_i, Z_i, S_i, Y_i)$ in a neighborhood of $\beta_0, \theta_0$. In addition, $E(M^2) < C_5 < \infty$ for all $n, d$.

(B4) The random variable $Z$ has compact support and its density $f_Z(z)$ is positive on that support. The bandwidth $h$ satisfies $nh^4 \to 0$ and $nh^2 \to \infty$. $\theta(z)$ has a bounded second derivative.

**Proof of Theorem 3.** Denote

$$J = [E\{(\mathcal{L}_\beta + \mathcal{L}_\theta \theta_\beta)|_{\beta_0, \theta_0}\}]^{-1}, \qquad \phi^*_{\mathit{eff}}(\beta, \theta) = J\mathcal{L}(\beta, \theta) \quad \text{and}$$
$$q'_{\lambda_n}(\beta) = Jp'_{\lambda_n}(\beta).$$

Let $\alpha_n = n^{-1/2} + a_n$ and $\phi^*_{\mathit{eff},i}(\beta, \hat{\theta}) = \phi^*_{\mathit{eff}}\{W_i, Z_i, S_i, Y_i, \beta, \hat{\theta}(\beta)\}$. It will be shown that

$$n^{-1/2} \sum_{i=1}^{n} \phi^*_{\mathit{eff},i}(\beta, \hat{\theta}) - n^{1/2} q'_{\lambda_n}(\beta) = 0 \tag{A4}$$

has a solution $\hat{\beta}$ that satisfies $\|\hat{\beta} - \beta_0\| = \mathrm{O}_p(d_n^{1/2}\alpha_n)$.

Due to the usual local estimating equation expansion, we have

$$\hat{\theta}(z, \beta_0) - \theta_0(z)$$
$$= (h^2/2)\theta_0''(z) - n^{-1}\sum_{j=1}^{n} K_h(Z_j - z)\Omega^{-1}(z)\Psi_j(\beta_0, \theta_0)/f_Z(z) \qquad (A5)$$
$$+ \mathrm{o}_p(n^{-1/2}),$$

which implies that $\hat{\theta}(z, \beta_0) - \theta_0(z) = \mathrm{O}_p(h^2 + n^{-1/2}h^{-1/2})$. For any $\beta$ such that $\|\beta - \beta_0\| = C\sqrt{d_n}\alpha_n$ for some constant $C$, we obtain the expansion

$$n^{-1/2}\sum_{i=1}^{n} \phi_{\mathit{eff},i}^*\{\beta, \hat{\theta}(\beta)\}$$

$$= n^{-1/2}\sum_{i=1}^{n} \phi_{\mathit{eff},i}^*\{\beta_0, \hat{\theta}(\beta_0)\}$$

$$+ n^{-1/2}\sum_{i=1}^{n} \left[\frac{\partial \phi_{\mathit{eff},i}^*\{\beta_0 + \hat{\theta}(\beta_0)\}}{\partial \beta^{\mathrm{T}}} + \frac{\partial \phi_{\mathit{eff},i}^*\{\beta_0 + \hat{\theta}(\beta_0)\}}{\partial \theta^{\mathrm{T}}}\frac{\partial \hat{\theta}}{\partial \beta^{\mathrm{T}}}\right](\beta - \beta_0)$$

$$+ \frac{1}{2\sqrt{n}}\sum_{i=1}^{n}(\hat{\beta} - \beta_0)^{\mathrm{T}}\frac{\mathrm{d}[\phi_{\mathit{eff},i}^*\{\beta_0 + \hat{\theta}(\beta_0)\} + \phi_{\mathit{eff},i}^*\{\beta_0 + \hat{\theta}(\beta_0)\}\frac{\mathrm{d}\hat{\theta}}{\mathrm{d}\beta}]}{\mathrm{d}\beta^{\mathrm{T}}}\bigg|_{\beta^*}(\beta - \beta_0),$$

where $\beta^*$ is between $\beta$ and $\beta_0$. Because of condition (B3), each component of the last term is uniformly of order $\mathrm{O}_p(n^{1/2}\|\beta - \beta_0\|^2)$. The second term can be written as $n^{1/2}\{1 + \mathrm{o}_p(1)\}(\beta - \beta_0)$ under conditions (B1), (B3) and (B4). The first term can be further expanded as

$$n^{-1/2}\sum_{i=1}^{n} \phi_{\mathit{eff},i}^*\{\beta_0, \hat{\theta}(\beta_0)\}$$

$$= n^{-1/2}\sum_{i=1}^{n} \frac{\partial \phi_{\mathit{eff},i}^*\{\beta_0, \hat{\theta}(\beta_0)\}}{\partial \theta^{\mathrm{T}}}\{\hat{\theta}(\beta_0) - \theta_0\} + \mathrm{O}_p(n^{1/2})\{\hat{\theta}(\beta_0) - \theta_0\}\{\hat{\theta}(\beta_0) - \theta_0\}^{\mathrm{T}}$$

$$= n^{-1/2}\sum_{i=1}^{n} \frac{\partial \phi_{\mathit{eff},i}^*\{\beta_0, \hat{\theta}(\beta_0)\}}{\partial \theta^{\mathrm{T}}}\{\hat{\theta}(\beta_0) - \theta_0\} + \mathrm{o}_p(1)$$

under conditions (B3) and (B4). Summarizing the above results, making use of (A5), we obtain

$$n^{-1/2}\sum_{i=1}^{n} \phi_{\mathit{eff},i}^*\{\beta, \hat{\theta}(\beta)\}$$

$$= n^{-1/2} \sum_{i=1}^{n} \phi^*_{\text{eff},i}(\beta_0, \theta_0) + n^{1/2}(\beta - \beta_0)$$

$$- n^{-3/2} \sum_{j,i=1}^{n} \frac{\partial \phi^*_{\text{eff},i}(\beta_0, \theta_0)}{\partial \theta^{\mathrm{T}}} K_h(Z_j - Z_i) \frac{\Omega(Z_i)\Psi_j(\beta_0, \theta_0)}{f_Z(Z_i)} + o_p(1)$$

$$= n^{-1/2} \sum_{i=1}^{n} \phi^*_{\text{eff},i}(\beta_0, \theta_0) + n^{1/2}(\beta - \beta_0) - n^{-1/2} \sum_{i=1}^{n} J\mathcal{U}(Z_i)\Psi_j(\beta_0, \theta_0) + o_p(1)$$

under condition (B4). Similar to the situation in Theorem 1, under condition (P1), we further obtain

$$(\beta - \beta_0)^{\mathrm{T}} \left\{ n^{-1/2} \sum_{i=1}^{n} \phi^*_{\text{eff},i}(\beta, \hat{\theta}) - n^{1/2} q'_{\lambda_n}(\beta) \right\}$$

$$= (\beta - \beta_0)^{\mathrm{T}} \left\{ n^{-1/2} \sum_{i=1}^{n} \phi^*_{\text{eff},i}(\beta_0, \theta_0) - n^{1/2} q'_{\lambda_n}(\beta_0) - n^{-1/2} \sum_{i=1}^{n} J\mathcal{U}(Z_i)\Psi_i(\beta_0, \theta_0) \right\} \quad (A6)$$

$$+ n^{1/2} \|\beta - \beta_0\|^2 + o_p\{n^{1/2}\|\beta - \beta_0\|^2\}.$$

The first term in the above display is of order $O_p(Cn^{1/2}d_n\alpha_n^2)$ and the second term equals $C^2 n^{1/2} d_n \alpha_n^2$, which dominates the first term as long as $C$ is large enough. The last term is dominated by the first two terms. Thus, for any $\epsilon > 0$, as long as $C$ is large enough, the probability for the above display to be larger than zero is at least $1 - \epsilon$. From Brouwer's fixed point theorem we know that with a probability of at least $1 - \epsilon$, there exists at least one solution for (A4) in the region $\|\beta - \beta_0\| \leq Cd_n^{1/2}\alpha_n$. □

*Lemma for Theorem 4.*

**Lemma 2.** *If conditions in Theorem 4 hold, then for any given $\beta$ that satisfies $\|\beta - \beta_0\| = O_p(\sqrt{d/n})$, with probability tending to 1, any solution $(\beta_I^{\mathrm{T}}, \beta_{II}^{\mathrm{T}})^{\mathrm{T}}$ of (2.2) satisfies that $\beta_{II} = 0$.*

**Proof.** Denote the $k$th equation in $\sum_{i=1}^{n} \mathcal{L}_i\{\beta, \hat{\theta}(\beta)\}$ as $L_{nk}(\beta, \hat{\theta})$ and that in $\sum_{i=1}^{n} \mathcal{U}(Z_i) \times \Psi_i(\beta_0, \theta_0)$ as $G_{nk}(\beta_0, \theta_0)$, $k = d_1 + 1, \ldots, d_n$, then the expansion in Theorem 3 leads to

$$L_{nk}(\beta, \hat{\theta}) - np'_{\lambda_n}(\beta_k)$$
$$= L_{nk}(\beta_0, \theta_0) - G_{nk}(\beta_0, \theta_0)$$
$$+ n \sum_{j=1}^{d} (J^{-1})_{kj}(\beta_j - \beta_{j0}) - np'_{\lambda_n}(|\beta_k|)\operatorname{sign}(\beta_k) + o_p(\sqrt{nd_n}).$$

Similar to the derivation in Lemma 1, the first three terms of the above display are all of order $O_p(\sqrt{nd_n})$, hence we have

$$L_{nk}(\beta,\hat{\theta}) - np'_{\lambda_n}(\beta_k) = -\sqrt{nd_n}\{\sqrt{n/d_n}p'_{\lambda_n}(|\beta_k|)\operatorname{sign}(\beta_k) + O_p(1)\}.$$

Because of (2.6), the sign of $L_{nk}(\beta) - np'_{\lambda_n}(\beta_k)$ is decided by $\operatorname{sign}(\beta_k)$ completely. From the continuity of $L_{nk}(\beta) - np'_{\lambda_n}(\beta_k)$, we obtain that it is zero at $\beta_k = 0$ with a probability larger than any $1 - \epsilon$. □

**Proof of Theorem 4.** (i) immediately follows by Lemma 1. Denote the first $d_1$ equations in $\sum_{i=1}^n \mathcal{L}_i\{(\beta_I^T, 0^T)^T, \hat{\theta}\}$ as $L_n\{\beta_I, \hat{\theta}(\beta_I)\}$ and that in $\sum_{i=1}^n \Psi_i(\beta_0, \theta_0)\mathcal{U}(Z_i)$ as $G_n(\beta_{I0}, \theta_0)$. Note that the $d_1 \times d_1$ upper left block of $J^{-1}$ is the matrix $A$ defined in Theorem 4. Using the Taylor expansion for the penalized estimating function at $\beta = (\beta_I^T, 0)^T$, the first $d_1$ equations yield

$$0 = L_n\{\hat{\beta}_I, \hat{\theta}(\hat{\beta}_I)\} - np'_{\lambda_n, I}(\hat{\beta}_I)$$
$$= L_n(\beta_{I0}, \theta_0) - G_n(\beta_{I0}, \theta_0) + nA(\hat{\beta}_I - \beta_{I0}) - nb_n$$
$$\quad - n\{\Sigma_n + o_p(1)\}(\hat{\beta}_I - \beta_{I0}) + o_p(d_n^{1/2}n^{1/2})$$
$$= L_n(\beta_{I0}, \theta_0) - G_n(\beta_{I0}, \theta_0) + n(A - \Sigma_n)[\hat{\beta}_I - \beta_{I0} - (A - \Sigma_n)^{-1}b_n] + o_p(d_n^{1/2}n^{1/2}).$$

Using condition (B2), we have

$$n^{1/2}v^T B^{-1/2}\{(-A + \Sigma_n)(\hat{\beta}_I - \beta_{I0}) + b_n\}$$
$$= n^{-1/2}v^T B^{-1/2}\{L_n(\beta_{I0}, \theta_0) - G_n(\beta_{I0}, \theta_0)\} + o_p(v^T B^{-1/2})$$
$$= n^{-1/2}v^T B^{-1/2}\{L_n(\beta_{I0}, \theta_0) - G_n(\beta_{I0}, \theta_0)\} + o_p(1).$$

Let $Y_i = n^{-1/2}v^T B^{-1/2}\{\mathcal{L}_{nIi}(\beta_{I0}, \theta_0) - \mathcal{U}_{nI}(Z_i)\Psi_i(\beta_{I0}, \theta_0)\}, i = 1, \ldots, n$. It follows that for any $\epsilon > 0$,

$$\sum_{i=1}^n E\|Y_i\|^2 \mathbf{1}(\|Y_i\| > \epsilon) = nE\|Y_1\|^2 \mathbf{1}(\|Y_1\| > \epsilon) \leq n(E\|Y_1\|^4)^{1/2}\{P_r(\|Y_1\| > \epsilon)\}^{1/2}.$$

Using the Chebyshev inequality, we have $P_r(\|Y_1\| > \epsilon) = O(n^{-1})$ and $E(\|Y_1\|^4)$ is bounded by

$$n^{-2}\lambda_{\max}^2(B^{-1})E[\{\mathcal{L}_{nI1}(\beta_{I0}, \theta_0) - \mathcal{U}_{nI}(Z_1)\Psi_1(\beta_{I0}, \theta_0)\}^T$$
$$\times \{\mathcal{L}_{nI1}(\beta_{I0}, \theta_0) - \mathcal{U}_{nI}(Z_1)\Psi_1(\beta_{I0}, \theta_0)\}]^2,$$

which equals $n^{-2}\lambda_{\max}^2(B^{-1})E\|\{\mathcal{L}_{nI1}(\beta_{I0}, \theta_0) - \mathcal{U}_{nI}(Z_1)\Psi_1(\beta_{I0}, \theta_0)\}\|^4 = O(d_n^2 n^{-2})$ by condition (B2). Hence,

$$\sum_{i=1}^n E\|Y_i\|^2 \mathbf{1}(\|Y_i\| > \epsilon) = O(nd_n n^{-1} n^{-1/2}) = o(1).$$

On the other hand,

$$\sum_{i=1}^{n} \text{cov}(Y_i) = n \text{cov}[n^{-1/2} v^{\text{T}} B^{-1/2} \{ \mathcal{L}_{nI1}(\beta_{I0}, \theta_0) - \mathcal{U}_{nI}(Z_1) \Psi_1(\beta_{I0}, \theta_0) \}] = 1.$$

(ii) follows by the Lindeberg–Feller central limit theorem. □

## Acknowledgements

Ma's work was supported by a Swiss NSF grant and a US NSF Grant DMS-0906341. Li's research was supported by a NSF Grant DMS-0348869 and National Institute on Drug Abuse (NIDA) Grants P50 DA-10075 and R21 DA-024260. The content is solely the responsibility of the authors and does not necessarily represent the official views of the NIDA or the National Institute of Health.

## References

Bickel, P.J. and Ritov, A.J.C. (1987). Efficient estimation in the errors-in-variables model. *Ann. Statist.* **15** 513–540. MR0888423

Cai, J., Fan, J., Li, R. and Zhou, H. (2005). Variable selection for multivariate failure time data. *Biometrika* **92** 303–316. MR2201361

Candès, E. and Tao, T. (2007). The Dantzig selector: Statistical estimation when $p$ is much larger than $n$ (with discussion). *Ann. Statist.* **35** 2313–2392. MR2382644

Carroll, R.J. and Hall, P. (1988). Optimal rates of convergence for deconvolving a density. *J. Amer. Statist. Assoc.* **83** 1184–1186. MR0997599

Carroll, R.J., Ruppert, D., Stefanski, L.A. and Crainiceanu, C. (2006). *Measurement Error in Nonlinear Models: A Modern Perspective*, 2nd ed. London: CRC Press. MR2243417

Delaigle, A. and Hall, P. (2007). Using SIMEX for smoothing-parameter choice in errors-in-variables problems. *J. Amer. Statist. Assoc.* **103** 280–287. MR2394636

Delaigle, A. and Meister, A. (2007). Nonparametric regression estimation in the heteroscedastic errors-in-variables problem. *J. Amer. Statist. Assoc.* **102** 1416–1426. MR2372541

Fan, J. (1991). On the optimal rates of convergence for nonparametric deconvolution problems. *Ann. Statist.* **19** 1257–1272. MR1126324

Fan, J. and Huang, T. (2005). Profile likelihood inferences on semiparametric varying-coefficient partially linear models. *Bernoulli* **11** 1031–1057. MR2189080

Fan, J. and Li, R. (2001). Variable selection via nonconcave penalized likelihood and its oracle properties. *J. Amer. Statist. Assoc.* **96** 1348–1360. MR1946581

Fan, J. and Lv, J. (2008). Sure independence screening for ultra-high dimensional feature space (with discussion). *J. Roy. Statist. Soc. Ser. B* **70** 849–911.

Fan, J. and Peng, H. (2004). Nonconcave penalized likelihood with a diverging number of parameters. *Ann. Statist.* **32** 928–961. MR2065194

Hall, P. and Ma, Y. (2007). Semiparametric estimators of functional measurement error models with unknown error. *J. Roy. Statist. Soc. Ser. B* **69** 429–446. MR2323761

Härdle, W., Liang, H. and Gao, J. (2000). *Partially Linear Models*. Heidelberg: Springer Physica. MR1787637

Hunter, D. and Li, R. (2005). Variable selection using MM algorithms. *Ann. Statist.* **33** 1617–1642. MR2166557

Kannel, W.B., Newton, J.D., Wentworth, D., Thomas, H.E., Stamler, J., Hulley, S.B. and Kjelsberg, M.O. (1986). Overall and coronary heart disease mortality rates in relation to major risk factors in 325,348 men screened for MRFIT. *Am. Heart J.* **112** 825–836.

Lam, C. and Fan, J. (2008). Profile-Kernel likelihood inference with diverging number of parameters. *Ann. Statist.* **36** 2232–2260. MR2458186

Li, R. and Liang, H. (2008). Variable selection in semiparametric regression modeling. *Ann. Statist.* **36** 261–286. MR2387971

Li, R. and Nie, L. (2007). A new estimation procedure for partially nonlinear model via a mixed effects approach. *Canad. J. Statist.* **35** 399–411. MR2396027

Li, R. and Nie, L. (2008). Efficient statistical inference procedures for partially nonlinear models and their applications. *Biometrics* **64** 904–911.

Liang, H., Härdle, W. and Carroll R.J. (1999). Estimation in a semiparametric partially linear errors-in-variables model. *Ann. Statist.* **27** 1519–1535. MR1742498

Liang, H. and Li, R. (2009). Variable selection for partially linear models with measurement errors. *J. Amer. Statist. Assoc.* **104** 234–248. MR2504375

Ma, Y. and Carroll, R.J. (2006). Locally efficient estimators for eemiparametric models with measurement error. *J. Amer. Statist. Assoc.* **101** 1465–1474. MR2279472

Ma, Y. and Li, R. (2007). Variable selection in measurement error models. Technical report. Available at http://www2.unine.ch/webdav/site/statistics/shared/documents/v10.pdf.

Ma, Y. and Tsiatis, A.A. (2006). Closed form semiparametric estimators for measurement error models. *Statist. Sinica* **16** 183–193. MR2256086

Severini, T.A. and Staniswalis, J.G. (1994). Quasilikelihood estimation in semiparametric models. *J. Amer. Statist. Assoc.* **89** 501–511. MR1294076

Stefanski, L.A. and Carroll, R.J. (1987). Conditional scores and optimal scores for generalized linear measurement-error models. *Biometrika* **74** 703–716. MR0919838

Tsiatis, A.A. and Ma, Y. (2004). Locally efficient semiparametric estimators for functional measurement error models. *Biometrika* **91** 835–848. MR2126036

Wang, H., Li, R. and Tsai, C. (2007). Tuning parameter selectors for the smoothly clipped absolute deviation method. *Biometrika* **94** 553–568. MR2410008

Zou, H. and Li, R. (2008). One-step sparse estimates in nonconcave penalized likelihood models (with discussion). *Ann. Statist.* **36** 1509–1566. MR2435443